\documentclass{article}  
\hoffset=,25cm 
\usepackage{parskip}
\usepackage{amsmath}
\usepackage[initials]{amsrefs}
\usepackage{amsfonts}
\usepackage{color}
\usepackage{amssymb,bbm}
\usepackage{eucal}
\usepackage{framed}
\usepackage{multicol}
\usepackage{graphicx}
\usepackage[makeroom]{cancel}
\usepackage{fancyhdr}
\usepackage{enumitem}
\usepackage{tikz}
\usepackage{hyperref}
\hypersetup{
  colorlinks   = true, 
  urlcolor     = blue, 
  linkcolor    = blue, 
  citecolor   = blue 
}
\usetikzlibrary{arrows}
\newcommand{\bc}{\begin{center}}
\newcommand{\ec}{\end{center}}
\newcommand{\beqn}{\begin{align*}}
\newcommand{\eeqn}{\end{align*}}
\newcommand{\benu}{\begin{enumerate}}
\newcommand{\eenu}{\end{enumerate}}
\newcommand{\bit}{\begin{itemize}}
\newcommand{\eit}{\end{itemize}}
\newcommand{\dsp}{\displaystyle}

\DeclareMathOperator{\tr}{tr}

\newcommand{\cfn}[1]{\ensuremath \raisebox{0pt}{$\mathbbm{1}$}_{#1}}

\newcommand{\en}{\rightarrow}

\newcommand{\bR}{\mathbb{R}}

\newcommand{\bZ}{\mathbb{Z}}

\newcommand{\bC}{\mathbb{C}}

\newcommand{\cR}{\mathcal{R}}

\newcommand{\cD}{\mathcal{D}}

\newcommand{\cM}{\mathcal{M}}

\newcommand{\eps}{\varepsilon}

\newcommand{\vect}[1]{{\vec{\mathbf{#1}}\mkern2mu}}

\newcommand{\limi}[3]{ {\displaystyle \lim_{#1 \rightarrow #2} #3} }

\newcommand{\subc}{\subseteq}

\newcommand{\pdd}[3]{\ifx#2#3\frac{\partial^2 #1}{\partial #2^2}\else \frac{\partial^2 #1}{\partial #2\,\partial #3}\fi}

\newcommand{\set}[2]{\left\{\,#1 : #2 \, \right\}}
\newcommand{\row}[3]{{#1}_{#2},\dots , {#1}_{#3}}

\newcommand{\inn}[2]{\left\langle #1 , #2 \right\rangle}
\newcommand{\com}[2]{\left[ #1, #2  \right]}
\newcommand{\hcb}[1]{\widehat{B}{(#1)}}
\newcommand{\hcf}[1]{\widehat{f}{(#1)}}

\newcommand{\lpn}[2]{\left\| #1 \right\|_{L^{#2}}}

\DeclareFontFamily{U}{wncy}{}
    \DeclareFontShape{U}{wncy}{m}{n}{<->wncyr10}{}
    \DeclareSymbolFont{mcy}{U}{wncy}{m}{n}
    \DeclareMathSymbol{\Sh}{\mathord}{mcy}{"58} 
\newcommand{\pa}[1]{\tilde{#1}}
\DeclareMathOperator{\Tr}{Tr}

\allowdisplaybreaks
\newtheorem{teo}{Theorem}[section]

\newtheorem{prop}{Proposition}[section]
\newtheorem{defi}{Definition}[section]
\newtheorem{lem}[teo]{Lemma}
\newcommand{\ip}[2]{\left\langle #1 , #2  \right\rangle}
\newcommand{\pp}[6]{ \sum_{R \in \cD^2} \ip{B}{h_R^{(#1,#2)}} \ip{f}{h_R^{(#3,#4)}} h_R^{(#5, #6)}}

\newcommand{\shc}[1]{ \com{\com{#1}{\Sh_1}}{\Sh_2} }

\begin{document}

\title{Characterization of two parameter matrix-valued BMO by commutator
with the Hilbert transform}

\author{Dar\'io Mena}
\date{}

\maketitle

\begin{abstract}
In this paper we prove that the space of two parameter, matrix-valued BMO functions can be
characterized by considering iterated commutators with the Hilbert transform.  Specifically, we
prove that
$$
\| B \|_{BMO} \lesssim \| [[M_B, H_1],H_2] \|_{L^2(\bR^2;\bC^d) \en L^2(\bR^2;\bC^d)} \lesssim \| B
\|_{BMO}.
$$
The upper estimate relies on Petermichl's
 representation of the Hilbert transform as an average of dyadic shifts, and the boundedness of
certain paraproduct operators, while the lower bound follows Ferguson and Lacey's proof for the
scalar case.
\end{abstract}

\section{Introduction}

It is well known, by the work of R. Coifman, R. Rochberg, and G. Weiss \cite{CRW1976}, that the space of functions of bounded mean oscillation (BMO) can
be
characterized by commutators with the Hilbert transform (and in general, with the Riesz transforms).  Given $b \in BMO$, let $M_b$ represent the
multiplication operator $M_b(f) = bf$, if $H$ represents the Hilbert transform, defined as
$$
H f(x) = p.v \frac{1}{\pi} \int_{\bR} \frac{f(y)}{x-y}\, dy,
$$
then we have
$$
\| b \|_{BMO} \lesssim \| [M_b, H] \|_{L^2 \en L^2} \lesssim \| b \|_{BMO}.
$$
The study of the norm of the commutator has several implications in the characterization of Hankel
operators, the problem of factorization and weak factorization of function spaces and the div-curl
problem.    Several extensions and generalizations have been made in different settings.  In the two
parameter version of this result, the upper bound was shown by S. Ferguson and C. Sadosky in
\cite{FeSa2000}, while the lower bound was proved by S. Ferguson
and M. Lacey in \cite{FeLa2002}.    
The formulation in this case is the following: If $H_i$ represents the Hilbert transform in the $i$-th variable, then
$$
\| b \|_{BMO} \lesssim \| [[M_b, H_1],H_2]\|_{L^2 \en L^2} \lesssim \| b \|_{BMO}.
$$
Here, we are considering the product BMO of S.Y. Chang and R. Fefferman \cite{ChFe1980}.  These results were later extended to the multi-parameter case by
M. Lacey and E. Terwilleger \cite{LaTe2009}.

The idea of the present work, is to obtain the same characterization in the two parameter case, for a matrix-valued BMO function.   In the one parameter
setting, we have the desired  characterization due to S. Petermichl \cite{Pe2000}, and also F. Nazarov, G. Pisier, S. Treil and A. Volberg \cite{NPTV2002}.

Consider the collection $\cD$ of dyadic intervals, that is 
$$
\cD : = \set{ [k 2^{-j}, (k+1)2^{-j}) }{j,k \in \bZ},
$$
and the collection of ``shifted'' dyadic intervals  
$$ 
\cD^{\alpha, r} = \set{\alpha + r[k2^j,
(k+1)2^j)}{k,j \in \bZ}, \quad \alpha, r \in \bR.
$$ 

Define the dyadic Haar function as $h_I : = \frac{1}{\sqrt{|I|}} (\cfn{I_-} - \cfn{I_+})$, where $I_-$ and $I_+$ represent the left and right half of the
interval $I$, respectively.   Denote also $h^1_J = \frac{\cfn{I}}{\sqrt{|I|}}$ (non-cancellative Haar function).     The family $\set{h_I}{I \in \cD}$ (or
$I \in \cD^{\alpha, r}$), is an orthonormal basis for $L^2(\bR; \bC^d)$; here, for two Banach spaces $X$ and $Y$, we use the notation $L^p(X;Y)$ to
denote the set $\set{f: X \en Y}{  \int_X \| f \|_Y^p < \infty}$.

Define the dyadic Haar shift by $\Sh^{\alpha, r} (h_I) =
\frac{1}{\sqrt{2}} (h_{I_-} - h_{I_+})$, and extend to a general function $f$ by 
$$\Sh^{\alpha, r}(f) = \sum_{I \in \cD^{\alpha, r}} \inn{f}{h_I} \Sh^{\alpha, r}(h_I) = \sum_{I \in \cD^{\alpha, r}} \inn{f}{h_I} \frac{1}{\sqrt{2}} (h_{I_-} - h_{I_+}).$$  

Note that $\Sh^{\alpha, r}$ is bounded
from $L^2(\bR; \bC^d)$ to $L^2(\bR; \bC^d)$, with operator norm $1$.
As proven by Petermichl in \cite{Pe2000}, the kernel for the Hilbert transform can be written as an average of dyadic shifts, in particular
$$
K(t, x) = \limi{L}{\infty}{  \frac{1}{2 \log L}  \int_{1/L}^L  \limi{R}{\infty}{  \frac{1}{2R} \int_{-R}^R K^{\alpha, r}(t, x) \, d\alpha } \,
\frac{dr}{r} }.
$$  

Where $K^{\alpha, r}(t, x) = \sum_{I \in \cD^{\alpha, r}} h_I(t) \Sh^{\alpha, r}(h_I (x))$.  Therefore, it is enough to prove the upper
bound for the commutator with the shift $[M_B, \Sh]$ (the estimates don't depend on $\alpha$ or $r$).  

Let $B$ be a function with values in the space of $d \times d$ matrices.  We consider the commutator $[M_B, H]$ acting on a vector-valued function $f$ by
$$
[M_B, H] f = B H(f) - H(Bf).
$$

The result obtained by Petermichl is based on a decomposition in paraproducts, and uses the estimates obtained by Katz \cite{Ka1997}, and Nazarov, Treil
and
Volberg \cite{NTV1997}
independently.  We have
$$
\left \|  [M_B , H]  \right \|_{L^2(\bR; \bC^d) \en L^2(\bR; \bC^d)}  \lesssim \log(1+d) \| B \|.
$$

Motivated by this result, we wish to find a generalization in a two parameter setting, with the corresponding definition of the product BMO space
(analogous to the one given by Chang and Fefferman in \cite{ChFe1980}). The main result of the paper can be stated as follows.

\begin{teo}
Let $B$ be a $d \times d$ matrix-valued BMO function on $\bR^2$.  If $M_B$ denotes the operator ``multiplication by $B$'', and $H_i$ represents the
Hilbert 
transform in the $i$-th parameter, for $i = 1, 2$, then the norm of the iterated commutator $[[M_B, H_1],H_2]$ satisfies
$$
d^{-2} \| B \|_{BMO} \lesssim \| [[M_B, H_1],H_2] \|_{L^2(\bR^2,\bC^d) \en L^2(\bR^2,\bC^d)}
\lesssim  d^3 \| B \|_{BMO}.
$$
\end{teo}
The paper is organized as follows. Section 2, contains the proof of the upper bound for the norm of
the commutator,  using a decomposition in paraproducts.  
Section 3 contains the proof of the lower bound, that relies on the proof for the scalar
case by S. Ferguson and M. Lacey in \cite{FeLa2002}.  Throughout the paper, we use the
notation $A \lesssim B$ to indicate that there is a
positive constant $C$, such that $A \leq C B$.

\bc
\textbf{Acknowledgment}  
\ec

The author wishes to thank his advisor Brett Wick for proposing this problem, and for all the fruitful conversations, his
many suggestions and guidance during the process.

\section{Upper bound}

Consider $\cR = \cD \times \cD$, the class of rectangles consisting on products of dyadic intervals.    Given a subset $E$ of $\bR^2$, denote by $\cR(E)$
the family of dyadic rectangles contained in $E$.

Consider the wavelet $w_I$ constructed by Meyer in \cite{Mey2000}, and the two-parameters wavelet $v_R(x,y) = w_I(x) w_J(y)$ for $R = I \times J$,  with
all its properties listed in \cite{FeLa2002}.   
We start by giving the definitions of product $BMO$ and product dyadic $BMO$.
\begin{defi}[BMO]\label{def:BMO}
A function $B$ is in $BMO(\bR^2)$ if and only if there are constants $C_1$ and $C_2$ such that, for
any open set $U \subc \bR^2$ we have 
\benu
\item $ \dsp  \left(   \frac{1}{|U|} \sum_{R \in \cR(U)} \inn{B}{v_R}\inn{B}{v_R}^*    \right)^{1/2}
\leq C_1 I_d $
\item $ \dsp \left(   \frac{1}{|U|} \sum_{R \in \cR(U)} \inn{B}{v_R}^*\inn{B}{v_R} 
\right)^{1/2} \leq C_2 I_d. $
\eenu
The inequalities are considered in the sense of operators, $I_d$ is the identity $d\times d$
matrix. The $BMO$-norm is defined as the smallest constant, denoted
by $\| B \|_{BMO}$, for which the two inequalities are satisfied simultaneously.  If we take the supremum only over 
rectangles $U$, we obtain the rectangular $BMO$-norm, denoted by $\| B \|_{BMO_{\textup{rec}}}$.  
\end{defi}

If $h_I$ represents the Haar function associated to a dyadic interval $I$, define
$$
h_R(x,y) = h_I(x)h_J(y), \quad \mbox{ for $R = I \times J$}.
$$
That is $h_R = h_I \otimes h_J$.  The family $\{ h_R \}_{R \in \cR}$ is an orthonormal basis for $L^{2}(\bR^2, \bC^d)$.  We have the following definition of dyadic $BMO$.  Note that it is the same definition, but considering the Haar wavelet instead of the Meyer wavelet.

\begin{defi}[Dyadic BMO]
A matrix-valued function $B$ is in $BMO_d(\bR^2)$ (dyadic $BMO$) if and only if, there are constants $C_1$ and $C_2$ such that for any open subset
$U$ of the plane, we have 
\benu
\item $\dsp \left ( \frac{1}{|U|} \sum_{R \in \cR(U)} \inn{B}{h_R} \inn{B}{h_R}^* \right )^{1/2}
\leq C_1
I_d $
\item $\dsp \left ( \frac{1}{|U|} \sum_{R \in \cR(U)} \inn{B}{h_R}^* \inn{B}{h_R}  \right
)^{1/2} \leq C_2 I_d.$
\eenu
Where the inequality is in the sense of operators. And the corresponding norm $\| B \|_{BMO_d}$ is,
again, the best constant for the two inequalities.
\end{defi}

It is known that $\| B \|_{BMO_d}  \leq \| B \|_{BMO}$; this fact can be found in \cite{Trei2008}.
In that paper, the proof of the inequality is given in the multiparameter setting, for Hilbert
space-valued functions, by means of the dual inequality $\| f \|_{H^1}  \leq \| f \|_{H^1_d}$
(Estimate 2.3 in \cite{Trei2008}).  The duality in the dyadic case is discussed later, in the proof
of Proposition \ref{paraproducts}.  Using this fact, for the
proof of the upper bound, it's enough to consider the dyadic version of $BMO$ for the
computations.   For the rest of this section, we use $\widehat{B}(R)$ to denote the Haar coefficient
of the function $B$, associated to the function $h_R$,
that is
$$
 \widehat{f}(R) = \inn{f}{h_R} = \int_{\bR^2} f(x,y) h_R(x,y) \, dx \, dy.
$$

Since $\widehat{B}(R) \widehat{B}(R)^*$ is a positive semi-definite matrix, we have
\begin{align*}
\sqrt{\frac{1}{|U|} \sum_{R \in \cR(U) } \| \widehat{B}(R)  \|^2} & \simeq \sqrt{ \Tr \left(
\frac{1}{|U|}\sum_{R \in \cR(U) }    \widehat{B}(R)
\widehat{B}(R)^* \right)}  \\
& \leq \Tr \sqrt{   \frac{1}{|U|}\sum_{R \in \cR(U) }    \widehat{B}(R)
\widehat{B}(R)^* }.	
\end{align*}

So, if we consider the two inequalities
$$ 
\sqrt{\frac{1}{|U|} \sum_{R \in \cR(U) } \widehat{B}(R) \widehat{B}(R)^*} \leq C  I_d,
\sqrt{\frac{1}{|U|} \sum_{R \in \cR(U) } \widehat{B}(R)^* \widehat{B}(R)} \leq C  I_d,
$$
taking the trace on both sides, we get
\begin{equation}\label{bmotrace}
\sqrt{\frac{1}{|U|} \sum_{R \in \cR(U) } \| \widehat{B}(R)  \|^2}  \leq  C d.  
\end{equation}

The initial computations are similar to the ones found in \cite{DaOu2014}.  In this, we need
simplified versions, since we are dealing only with the biparameter Hilbert transform; differences
will arise when we deal with the various paraproducts that result from this process, due to the
$BMO$ symbol being a matrix (which implies losing commutativity and requiring the use of matrix
norms). Similar computations are used in \cite{LPPW2010}, and this ideas can also be implemented in
our case. Although we can use some equivalent results from \cite{MPTT2004,MPTT2006} to deal with the
boudnedness of the paraproducts, the ones arising from our computations can be given self
contained proofs of their boundedness.

The dyadic shift operator $\Sh(f) = \sum_{I \in
\cD} \widehat{f}(I) \frac{1}{\sqrt{2}}(h_{I_{-}} -
h_{I_+})$  corresponds to the operator $S^{1,0}$ described by Dalenc and Ou in \cite{DaOu2014},
given by
$$
S^{1,0} f = \sum_{K \in \cD} \sum_{I \subc K}^{(0)} \sum_{J \subc K}^{\tiny (1)} a_{IJK} \langle f,
h_I \rangle h_J, \quad
a_{IJK} = \begin{cases}
          \frac{1}{\sqrt{2}} , & \mbox{if $J = K_-$}, \\
	  -\frac{1}{\sqrt{2}}, & \mbox{if $J = K_+$}.
          \end{cases}
$$
Here, the symbol $\sum_{I \subc J}^{(k)}$ represents summing over those dyadic intervals $I$ such that $I \subc J$, and $|I| =
2^{-k}|J|$.
Let $\tilde{I}$ represent the parent of the dyadic interval $I$, that is, the unique dyadic interval containing $I$ with $|\tilde{I}| = 2|I|$,  then, the shift can also be expressed in a simpler way by
\begin{equation}\label{shift}
\Sh(f) = \sum_{I \in \cD} a_I \widehat{f}(\tilde{I}) h_I,
\end{equation}
where $a_I = \frac{1}{\sqrt{2}}$ if $I = \tilde{I}_{-}$, and $- \frac{1}{\sqrt{2}}$ if $I = \tilde{I}_{+}$.

If we write $B = \sum_{I \in \cD} \widehat{B}(I) h_I$, and $f = \sum_{J \in \cD}
\widehat{f}(J) h_J$, then we can write 
$$
Bf = \sum_{I} \sum_{J} \widehat{B}(I)h_I \widehat{f}(J)h_J.
$$
Therefore the commutator 
$$
[M_B, \Sh](f) = M_B  \Sh (f) - \Sh(M_B f) = B \Sh(f) - \Sh( B f),
$$
can be written as
\begin{align*}
[M_B, \Sh](f) & = \sum_{I,J} \widehat{B}(I) \widehat{f}(J)h_I\Sh(h_J) - \sum_{I,J} \widehat{B}(I)
\widehat{f}(J)\Sh(h_I h_J) \\
& = \sum_{I,J} \widehat{B}(I)
\widehat{f}(J)[M_{h_I}, \Sh](h_J).	
\end{align*}

Note that the terms are non-zero, only when $I \cap J \neq \emptyset$, also, if $J \subsetneq I$, we
have that $h_I$ is constant in $I \cap J$, therefore,
for every $x \in I \cap J$, we have
\begin{align*}
[M_{h_I}, \Sh](h_J) & = h_I(x)\Sh(h_J(x)) - \Sh(h_I(x) h_J(x)) \\
& = h_I(x)\Sh(h_J(x)) - h_I(x)\Sh(h_J(x)) = 0.
\end{align*}
Then, the only non-trivial terms are those for which $I \subset J$.

We consider the two parameter commutator $[[M_B, H_1],H_2]$ acting on a vector-valued function $f$ by
\begin{align*}
[[M_B, H_1],H_2] f  = & B H_1 (H_2(f)) - H_1(B (H_2(f)) ) \\ 
& - H_2(B H_1(f)) + H_2 (H_1(Bf)).
\end{align*}
Where $H_1$ and $H_2$ represent the Hilbert transform, on the first and second variable respectively.  That is,
$$
H_1 f(x , y) = p.v \frac{1}{\pi} \int_{\bR} \frac{f(z,y)}{x-z}\, dz, \qquad  H_2 f(x , y) = p.v \frac{1}{\pi} \int_{\bR} \frac{f(x,z)}{y-z}\, dz.
$$
The main result we want to prove in this section is the following
\begin{teo}\label{upperbound}
Let $B$ be a matrix-valued $BMO_d(\bR^2)$ function and $f$ in $L^2(\bR^2; \bC^d)$, then
$$
\| [[M_B, H_1], H_2] \|_{L^2(\bR^2; \bC^d) \en L^2(\bR^2; \bC^d)} \lesssim \| B \|_{BMO_d}.
$$
\end{teo}

{\bf Proof: } Let $\Sh_1$ and $\Sh_2$ represent the dyadic shift operator in the first and second variable respectively, that is, $\Sh_1(h_R) = \Sh(h_I)
\otimes h_J $, and $\Sh_2(h_R) = h_I \otimes \Sh(h_J)$, for $R =  I \times J$, and extending to a function $f$ by
$$
\Sh_j(f) = \sum_{R \in \cR} \hcf{R} \Sh_j(h_R), \quad j = 1,2.
$$  
Or in the notation of (\ref{shift}),
$$
\Sh_1(f) = \sum_{I,J \in \cD} a_I \widehat{f}(\tilde{I}\times J) h_I \otimes h_J, \qquad  \Sh_2(f) = \sum_{I,J \in \cD} a_J \widehat{f}(I \times \tilde{J}) h_I \otimes h_J.
$$

Again, due to the representation of $H$ as an average of shifts, it is enough to prove the result for the commutator $\com{\com{M_B}{\Sh_1}}{\Sh_2}$.   By
an iteration of the computation for the one parameter case, using the Haar expansion of the functions $B$ and $f$ and taking their formal product, we obtain that $\com{M_B}{\Sh_1} (f)$ is equal to
\begin{align*}
& \sum_{R,S \in \cR} \widehat{B}(R) \widehat{f}(S) \left( h_R \Sh_1 (h_S) - \Sh_1(h_R h_S) \right) \\ 
= & \sum_{R,S \in \cR}
\widehat{B}(R) \widehat{f}(S) \com{M_{h_R}}{\Sh_1}(h_S) \\
= &  \sum_{I,J,K,L \in \cD} \widehat{B}(I \times J) \widehat{f}(K \times L)\left (  h_I \Sh_1 h_K   - \Sh_1 (h_I h_K)  \right )\otimes
h_J h_L.
\end{align*}
Repeating the same computations, we get that the two parameters commutator $\com{\com{M_B}{\Sh_1}}{\Sh_2} (f)$ is equal to
\begin{align*}
&   \sum_{I,J \in \cD} \sum_{K,L \in \cD} \widehat{B}(I \times J) \widehat{f}(K \times L)  h_I \Sh_1 h_K \otimes h_J
\Sh_2 h_L   \\
 - &  \sum_{I,J \in \cD} \sum_{K,L \in \cD} \widehat{B}(I \times J) \widehat{f}(K \times L) \Sh_1 (h_I h_K) \otimes h_J \Sh_2 h_L  \\
- &  \sum_{I,J \in \cD} \sum_{K,L \in \cD} \widehat{B}(I \times J) \widehat{f}(K \times L)h_I \Sh_1 h_K \otimes  \Sh_2 (h_Jh_L) \\
+ &   \sum_{I,J \in \cD} \sum_{K,L \in \cD} \widehat{B}(I \times J) \widehat{f}(K \times L)\Sh_1 (h_I h_K) \otimes \Sh_2 (h_Jh_L) \\
= &  T_1f - T_2f - T_3f + T_4f \\
= &  \sum_{I,J \in \cD} \sum_{K,L \in \cD} \widehat{B}(I \times J) \widehat{f}(K \times L) \com{M_{h_I}}{\Sh_1}(h_K) \otimes \com{M_{h_J}}{\Sh_2}(h_L).
\end{align*}
If either $I \cap K = \emptyset $, $J \cap L = \emptyset$, $K \subsetneq I$ or $L \subsetneq J$,  then we  have that $\com{M_{h_I}}{\Sh_1}(h_K) \otimes
\com{M_{h_J}}{\Sh_2}(h_L) = 0$; therefore, the terms are non-trivial only when $I \subc K$ and $J \subc L$.   We have four different cases, that can be
analyzed independently for each term in the sum. The computations for the four terms are similar, only the complete details for the term $T_2$ will be
provided, and at the end of the proof of the proposition we mention briefly how to deal with the other cases.    Let $\tilde T_j$ represent $T_j$ restricted to the case $I \subc K$ and $J \subc L$, then we have.
$$
\tilde T_2 f  = \Sh_1 \left ( \sum_{K} \sum_{L} \sum_{I \subc K} \sum_{J \subc L} \widehat{B}(I \times J) \widehat{f}(K \times L) h_I h_K \otimes h_J \Sh_2 h_L \right ).
$$
To analyze each of the four cases, we need the following proposition.

\begin{prop}\label{paraproducts}
Consider the following paraproducts
\benu
\item[(i)] $\dsp P^1_B(f) = \sum_{I,J \in \cD} \pm \widehat{B}(I \times \tilde{J})  \inn{f}{h_I\otimes h_{J}} h_I^1 \otimes  h_{J}
|I|^{-1/2} |\tilde{J}|^{-1/2}$.
\item[(ii)] $\dsp P^2_B(f) = \sum_{I,J}  \pm \widehat{B}(I \times \tilde{J})  \inn{f}{h_I^1 \otimes h_{\tilde{J}}}  h_I  \otimes  h_J
|I|^{-1/2} |\tilde{J}|^{-1/2}$.
\item[(iii)] $\dsp P^3_B(f) = \sum_{I,J \in \cD}  \widehat{B}(I \times J) \inn{ f}{h_I^1 \otimes h_J^1} h_I  \otimes  h_J |I|^{-1/2}
|J|^{-1/2}$. 
\item[(iv)] $\dsp P^4_B(f) = \sum_{I,J \in \cD}  \widehat{B}(I \times J)   \inn{ f }{h_I \otimes h^1_J} h_I^1 \otimes  h_J |I|^{-1/2}
|J|^{-1/2}$.
\item[(v)] $\dsp P^5_B(f) = \sum_{I,J \in \cD}  \widehat{B}(I \times J)   \inn{ f }{h^1_I \otimes h_J} h_I \otimes  h_J^1 |I|^{-1/2}
|J|^{-1/2}$.
\eenu
We have that for $i = 1,2,3,4$, 
$$\| P^i_B(f) \|_{L^2(\bR^2; \bC^d)} \lesssim d \| B \|_{BMO_d} \| f \|_{L^2(\bR^2; \bC^d)}.$$
\end{prop}

{\bf Proof of proposition: }  In the following computations, for simplification we will write
$L^2(Y) = L^2(\bR^2;Y)$, since all the functions that we consider are defined on $\bR^2$.

\textbf{(i)} We make use of a well known result, which is discussed in \cite{Ch1979} for the bidisc
case, but it is easily extended to the plane.

\begin{teo}[Carleson Embedding Theorem]\label{CET}
Let $\{ a_R \}_{R \in \cR}$ be a sequence of nonnegative numbers, indexed by the grid of dyadic rectangles.  Then the following are
equivalent:
\benu
\item [(i)] $\sum_{R \in \cR} a_R \langle f \rangle_R^2  \leq C_1 \lpn{f}{2}^2$, for all $f \in  L^2$.
\item [(ii)] $ \frac{1}{|U|} \sum_{R \in \cR(U)} a_R \leq C_2$, for all connected open sets $U
\subc \bR^2$.
\eenu
Moreover, $C_1 \simeq C_2$.
\end{teo}

We have the following basic estimates
\begin{align*}
& \left | \inn{P^1_B f}{g}_{L^2} \right |  = \left | \int_{\bR^2}  \inn{P^1_B f}{g}_{\bC^d} \, dx \, dy \right | \\
= &  \left | \int_{\bR^2}  \inn{\sum_{I,J} \pm \widehat{B}(I \times \pa{J})  \widehat{f}(I \times J)  \cfn{I} |I|^{-1} \otimes 
h_{J}
|\tilde{J}|^{-1/2}}{g}_{\bC^d} \, dx \, dy \right | \\
= &  \left | \int_{\bR^2} \sum_{I,J} \inn{\pm \widehat{B}(I \times \pa{J})  \widehat{f}(I \times J)  }{g  \cfn{I} |I|^{-1} \otimes 
h_{J}
|\tilde{J}|^{-1/2}}_{\bC^d} \, dx \, dy \right | \\
= &  \left | \sum_{I,J} \int_{\bR^2}  \inn{\pm \widehat{B}(I \times \pa{J})  \widehat{f}(I \times J)  }{g  \cfn{I} |I|^{-1} \otimes 
h_{J}
|\tilde{J}|^{-1/2}}_{\bC^d} \, dx \, dy \right | \\
= &  \left | \sum_{I,J}  \inn{\pm \widehat{B}(I \times \pa{J})  \widehat{f}(I \times J)  }{\int_{\bR^2}  \frac{1}{\sqrt{2}}g 
\cfn{I}
|I|^{-1} \otimes  h_{J} |J|^{-1/2}  \, dx \, dy}_{\bC^d}  \right | \\
= &   \frac{1}{\sqrt{2}} \left | \sum_{I,J}  \inn{\pm \widehat{B}(I \times \pa{J})  \widehat{f}(I \times J)  }{\inn{g} {  \cfn{I}
|I|^{-1} \otimes 
h_{J} |J|^{-1/2} }}_{\bC^d}  \right | \\
\leq &   \frac{1}{\sqrt{2}}  \sum_{I,J}  \left | \inn{\pm \widehat{B}(I \times \pa{J})  \widehat{f}(I \times J)  }{\inn{g} { 
\cfn{I} |I|^{-1}
\otimes  h_{J} |J|^{-1/2} }}_{\bC^d}  \right | \\
\leq &  \frac{1}{\sqrt{2}}   \sum_{I,J}  \left \|  \widehat{B}(I \times \pa{J}) \right \|  \left \| \widehat{f}(I \times J) \right
\|_{\bC^d} \left
\|\inn{g} {  \cfn{I} |I|^{-1} \otimes  h_{J} |J|^{-1/2} }\right \|_{\bC^d}  \\
\leq &  \frac{1}{\sqrt{2}}   \sum_{I,J} \left \| \widehat{f}(I \times J) \right \|_{\bC^d} \left \|  \widehat{B}(I \times \pa{J})
\right \|  \left
\langle \left \| g  \right \|_{\bC^d} \right \rangle_{I \times J} \\
\leq &  \frac{1}{\sqrt{2}}   \left ( \sum_{I,J} \left \| \widehat{f}(I \times J)  \right \|^2_{\bC^d} \right )^{\frac{1}{2}}  \left
( \sum_{I} 
\sum_{J} \left \|  \widehat{B}(I \times \pa{J}) \right \|^2  \left \langle \left \| g  \right \|_{\bC^d} \right \rangle_{I \times J}^2
\right )^{\frac{1}{2}}\\
\leq &  \frac{1}{\sqrt{2}} \| f \|_{L^2(\bC^d)}  \left ( \sum_{I,J}  \left \|  \widehat{B}(I \times \pa{J}) \right \|^2   \left \langle
\left \| g  \right \|_{\bC^d} \right \rangle_{I \times J}^2 \right)^{\frac{1}{2}}  \\
\lesssim & \| f \|_{L^2(\bC^d)} d \| B \|_{BMO_d} \| \| g \|_{\bC^d} \|_{L^2(\bR)} =  d\| B \|_{BMO_d}\| f \|_{L^2(\bC^d)}\| g
\|_{L^2(\bC^d)}. \\
\end{align*}
Here, we used the fact that since $B \in BMO_d$, then by (\ref{bmotrace}), the second condition in Theorem \ref{CET} is satisfied with $a_R = 
\left\| \hcb{R} \right\|^2$.  Note, that we have a linear dependence on the dimension of the matrix, due to the use of the trace.  Note also that the same
computations allow us to replace each individual $I$ and $J$ for a parent or ``great parent'' of $I$
and $J$, in which case, the implied constant will
depend also on complexity (level of relation with its ancestor); we will use
$P^1_B$ to denote any of these kind of paraproducts.

\textbf{(ii)} A direct computation shows that $(P^2_B)^*$ is of the type $P^1_{B^*}$, therefore, by the symmetry of the definition of
$BMO_d$-norm, the boundedness for $P^2_B$ follows from that of $P^1_B$.

\textbf{(iii)}  Denote by $S_2^d$ the space of $d \times d$
complex matrices, equipped with the norm derived from the inner product $\inn{A}{B}_{\Tr} = \tr({AB^*})$, that is $\| A \|_{S_2^d}^2  = \tr(AA^*) $. To
estimate the $L^2$-norm of this operator, we compute $\inn{P^3_B( f )}{g}$.
\begin{align*}
 & = \int_{\bR^2} \inn{\sum_{I,J} \widehat{B}(I \times J)   \inn{ f }{h_I^1 \otimes h^1_J} \frac{h_I \otimes  h_J}{|I|^{\frac{1}{2}}
|J|^{\frac{1}{2}}}}{g}_{\bC^d} \,
dx \, dy \\
& =\sum_{I,J}  \int_{\bR^2}  \inn{ \widehat{B}(I \times J)   \inn{ f }{h_I^1 \otimes h^1_J} }{g  \frac{h_I \otimes  h_J}{|I|^{\frac{1}{2}}
|J|^{\frac{1}{2}}} }_{\bC^d} \,
dx \, dy \\ 
& =\sum_{I,J} \inn{ \widehat{B}(I \times J)   \inn{ f }{h_I^1 \otimes h^1_J} }{\inn{g}{h_I \otimes  h_J}\frac{1}{|I|^{\frac{1}{2}}
|J|^{\frac{1}{2}}} }_{\bC^d}  \\ 
& =\sum_{I,J} \inn{ \widehat{B}(I \times J) }{\inn{g}{h_I \otimes  h_J} \inn{ f }{h_I^1 \otimes h^1_J}^*\frac{1}{|I|^{\frac{1}{2}}
|J|^{\frac{1}{2}}} }_{\Tr}  \\ 
& =\sum_{I,J} \int_{\bR^2} \inn{ B h_I \otimes h_J }{\inn{g}{h_I \otimes  h_J} \inn{ f }{h_I^1 \otimes h^1_J}^*\frac{1}{|I|^{\frac{1}{2}}
|J|^{\frac{1}{2}}} }_{\Tr} \, dx\, dy 
\\ 
& =\int_{\bR^2} \inn{ B  }{ \sum_{I,J} \inn{g}{h_I \otimes  h_J} \inn{ f }{h_I^1 \otimes h^1_J}^* \frac{h_I \otimes  h_J}{|I|^{\frac{1}{2}}
|J|^{\frac{1}{2}}} }_{\Tr} \, dx\, dy 
\\ 
& = \inn{ B  }{ \sum_{I,J} \inn{g}{h_I \otimes  h_J} \inn{ f }{h_I^1 \otimes h^1_J}^* \frac{h_I \otimes  h_J}{|I|^{\frac{1}{2}}
|J|^{\frac{1}{2}}} }_{L^2(S_2^d)} \\
& = \inn{B}{\Pi_1(f,g)}.
\end{align*}

Define the space $H^1_{d}$ to be the space of $d\times d$ matrix-valued functions $\Phi$ such that
$\| \Phi \|_{H^1_{d}} = \| S \Phi \|_{L^1}$, where $S$ is the square function defined by
$$
S^2 \Phi  (x , y) : = \sum_{I \in \cD} \sum_{J \in \cD} \|  \inn{\Phi}{h_I \otimes h_J} 
\|^2_{S_2^d} \frac{\cfn{I}(x)}{|I| } \frac{\cfn{J}(y)}{|J|}.
$$
Note that if $\Phi$ is in $H^1_{d}$, then all of its components are in scalar $H^1$, and for 
$1 \leq i, j \leq d$, we have $\| \Phi_{i,j} \|_{H^1} \leq \| \Phi \|_{H^1_d}$.   Also, if
$B$ is a matrix-valued $BMO_d$ function, then all of its components are in scalar dyadic
$BMO$, and an easy computation shows that for $1 \leq i, j \leq d$, $\| B_{i,j} \|_{BMO} \leq d \| B
\|_{BMO_d}$. Using these facts, we can easily verify the following duality statement:

\begin{lem}[$\mathbf{BMO_d-H^1_{d}}$ duality]
Let $B$ in $BMO_d$ and $\Phi$ in $H^1_{d}$, then
$$
\inn{ B }{ \Phi }_{L^2(S_2^{d})} \lesssim d^3 \| B \|_{BMO_d} \| \Phi \|_{H^1_{d}}.
$$
\end{lem}

Using this result, it is enough to prove that 
$$\left\| \Pi_1(f,g) \right\|_{H^1_{d}} \simeq \left\| S(\Pi_1(f,g)) \right\|_{L^1} \lesssim \| f
\|_{L^2} \| g \|_{L^2}.$$
We compute $[S(\Pi_1(f,g))(x,y)]^2$ to get
\begin{align*}
 & = \sum_{I,J} \left\|  \inn{g}{h_I \otimes h_J}  \inn{f}{h_I^1 \otimes
h_J^1}^* |I|^{-1/2} |J|^{-1/2}\right\|^2_{S_2^d} \frac{\cfn{I}(x)
\cfn{J}(y)}{|I||J|}
 \\
& =  \sum_{I,J} \left\|  \inn{g}{h_I \otimes h_J}  \right\|_{\bC^d}^2  \left\| \inn{f}{ \frac{h_I \otimes  h_J}{|I|^{\frac{1}{2}}
|J|^{\frac{1}{2}}} } \right\|^2_{\bC^d}
\frac{\cfn{I\times J}(x,y) }{|I \times J|} \\
& \leq  \sup_{(x,y) \in I\times J}  \left\| \inn{f}{ \frac{h_I \otimes  h_J}{|I|^{\frac{1}{2}}
|J|^{\frac{1}{2}}} } \right\|^2_{\bC^d}  \sum_{I,J} \left\|  \inn{g}{h_I \otimes
h_J}
\right\|_{\bC^d}^2 \frac{\cfn{I\times J}(x,y) }{|I \times J|} \\
& \leq  \sup_{(x,y) \in I\times J}   \inn{\left\| f\right\|_{\bC^d} }{ \frac{h_I \otimes  h_J}{|I|^{\frac{1}{2}}
|J|^{\frac{1}{2}}} } ^2  \sum_{I,J} \left\|  \inn{g}{h_I
\otimes h_J}
\right\|_{\bC^d}^2 \frac{\cfn{I\times J}(x,y) }{|I \times J|} \\
& \leq  \left[ \cM(\| f \|_{\bC^d}) (x,y)\right]^2  [S(g)(x,y)]^2.  \\
\end{align*}
Here, $\cM$ represents the strong maximal function.  Using the $L^2$-boundedness of the maximal and square functions, we conclude
$$
\left\| \Pi_1(f,g) \right\|_{H^1_{d}} \lesssim \left\| S(\Pi_1(f,g)) \right\|_{L^1} \lesssim 
\left\|  \cM(\| f \|_{\bC^d})  S(g)  \right\|_{L^1} \lesssim \| f
\|_{L^2} \| g \|_{L^2}.
$$

\textbf{(iv)}  As in the previous case, we compute $\inn{P^4_B( f )}{g}$
\begin{align*}
 & = \int_{\bR^2} \inn{\sum_{I,J} \widehat{B}(I \times J)   \inn{ f }{h_I \otimes h^1_J} \frac{h_I^1 \otimes  h_J}{|I|^{\frac{1}{2}} |J|^{\frac{1}{2}}}}{g}_{\bC^d} \,
dx \, dy \\
& =\sum_{I,J}  \int_{\bR^2}  \inn{ \widehat{B}(I \times J)   \inn{ f }{h_I \otimes h^1_J} }{g  \frac{h_I^1 \otimes  h_J}{|I|^{\frac{1}{2}} |J|^{\frac{1}{2}}} }_{\bC^d} \,
dx \, dy \\ 
& =\sum_{I,J} \inn{ \widehat{B}(I \times J)   \inn{ f }{h_I \otimes h^1_J} }{\inn{g}{h_I^1 \otimes  h_J}\frac{1}{|I|^{\frac{1}{2}} |J|^{\frac{1}{2}}} }_{\bC^d}  \\ 
& =\sum_{I,J} \inn{ \widehat{B}(I \times J) }{\inn{g}{h_I^1 \otimes  h_J} \inn{ f }{h_I \otimes h^1_J}^*\frac{1}{|I|^{\frac{1}{2}} |J|^{\frac{1}{2}}} }_{\Tr}  \\ 
& =\sum_{I,J} \int_{\bR^2} \inn{ B h_I \otimes h_J }{\inn{g}{h_I^1 \otimes  h_J} \inn{ f }{h_I \otimes h^1_J}^*\frac{1}{|I|^{\frac{1}{2}} |J|^{\frac{1}{2}}} }_{\Tr} \, dx\, dy 
\\ 
& =\int_{\bR^2} \inn{ B  }{ \sum_{I,J} \inn{g}{h_I^1 \otimes  h_J} \inn{ f }{h_I \otimes h^1_J}^* \frac{h_I \otimes  h_J}{|I|^{\frac{1}{2}} |J|^{\frac{1}{2}}} }_{\Tr} \, dx\, dy 
\\ 
& = \inn{ B  }{ \sum_{I,J} \inn{g}{h_I^1 \otimes  h_J} \inn{ f }{h_I \otimes h^1_J}^* \frac{h_I \otimes  h_J}{|I|^{\frac{1}{2}} |J|^{\frac{1}{2}}} }_{L^2(S_2^d)} \\
& = \inn{B}{\Pi_2(f,g)}. 
\end{align*}
Therefore, by duality, it is enough to prove that 
$$\left\| \Pi_2(f,g) \right\|_{H^1_{d}}  \lesssim \| f \|_{L^2} \| g \|_{L^2}.$$

For this, we proceed again to find a pointwise estimate for the square function.  We compute $[S(\Pi_2(f,g))]^2$
\begin{align*}
& = \sum_{I,J} \left\|   \inn{g}{h_I^1 \otimes h_J}  \inn{f}{h_I \otimes h_J^1}^*  \frac{1}{|I|^{\frac{1}{2}} |J|^{\frac{1}{2}}}\right\|^2_{S_2^d} \frac{\cfn{I
\times J}}{|I \times J|}  \\
& =  \sum_{I,J} \left\|   \left\langle \inn{g}{ h_J} \right\rangle_I   \right\|^2_{\bC^d} \frac{\cfn{J}}{|J|}  \left\|   \left\langle \inn{f}{ h_I}
\right\rangle_J   \right\|^2_{\bC^d} \frac{\cfn{I}}{|I|}   \\
& \leq \sum_{I,J}    \left\langle \left\| \inn{g}{ h_J} \right\|_{\bC^d} \right\rangle_I^2 \frac{\cfn{J}}{|J|}    \left\langle \left\| \inn{f}{ h_I}
\right\|_{\bC^d} \right\rangle_J^2 \frac{\cfn{I}}{|I|}   \\ 
& \leq \left(  \sum_{I}   \left( \cM_2 \left\| \inn{f}{ h_I}\right\|_{\bC^d} \right)^2 \frac{\cfn{I}}{|I|} \right)  \left(  \sum_{J} \left(
\cM_1 \left\|
\inn{g}{ h_J} \right\|_{\bC^d} \right)^2 \frac{\cfn{J}}{|J|} \right).
\end{align*}
Where $\cM_1$ and $\cM_2$ represent the maximal function in the first and second variable, respectively.  
These last two factors are symmetric to each
other, so it is enough to prove the $L^2$-boundedness for the operator
$$
\widetilde{S}f(x,y) = \left( \sum_{I}   \left( \cM_2 \left\| \inn{f}{ h_I}\right\|_{\bC^d}(y) \right)^2 \frac{\cfn{I}(x)}{|I|} \right)^{1/2}.
$$
But this is easy, since
\begin{align*}
\int_{\bR^2} (  \widetilde{S}f(x,y)  )^2 \, dx \, dy  & =  \sum_{I}  \int_{\bR} \left( \cM_2
\left\| \inn{f}{ h_I}\right\|_{\bC^d}(y) \right)^2 \, dy  \\
& \lesssim  \sum_{I}  \int_{\bR} \left\|
\inn{f(\cdot, y)}{ h_I(\cdot)}\right\|_{\bC^d}^2 \, dy = \| f \|_{L^2}^2.
\end{align*}

\textbf{(v)} The computatios are symmetric to those for \textbf{(iv)}, exchanging the roles of $I$
and $J$.  \hfill $\blacksquare$

We proceed now to prove the upper bound for the four different cases.  In each of them, the idea is to reduce the term to an expression of the form
$\Sh_1 \circ P_B^i \circ \Sh_2$, therefore, by Proposition \ref{paraproducts} and the boundedness of the shifts, we get the desired result.  The estimates for the rest of the terms are similar, since they
are reduced to find an upper bound for the norm of the four variants of paraproduct studied above.  More specifically, they correspond to expressions of the form $\Sh_i (P_B(\Sh_j f))$, $\Sh_i (\Sh_j (P_Bf))$  and $\Sh_i (\Sh_j (P_Bf))$, $\Sh_i(P_B f)$, or duals of operators of the form $\Sh_i (P_{B^*}(\Sh_j f))$, $\Sh_i(\Sh_j (P_{B^*}f))$, $\Sh_i (\Sh_j (P_{B^*}f))$ and $\Sh_i(P_{B^*} f)$.

{\bf Case $\mathbf{I = K}$, $\mathbf{J = L}$.}  In this case, using the definition of the shift, we have
\begin{align*}
 & \Sh_1 \left (  \sum_{I} \sum_{J} \widehat{B}(I \times J) \widehat{f}(I \times J) h_I^2   h_J \Sh_2 h_J \right ) \\ 
=  & \Sh_1 \left (  \sum_{I}  \sum_{J}  \widehat{B}(I \times \tilde{J}) \widehat{f}(I \times \tilde{J}) h_I^2 \otimes {h_{\tilde{J}}} a_J h_{J} \right ).
\end{align*}

Since $\Sh_2 \inn{f}{h_I} = \sum_L a_L \widehat{f}(I \times \pa{L}) h_L$, then, $\inn{\Sh_2 \inn{f}{h_I}}{h_{J}} = a_J \widehat{f}(I \times \pa{J})$.  So, the previous expression is equal to
\begin{align*}
& \Sh_1 \left (    \sum_{I}  \sum_{J} \widehat{B}(I \times \tilde{J})  \inn{\Sh_2 \inn{f}{h_I}}{h_{J}} h_I^2 \otimes h_{\tilde{J}} h_{J} \right ) \\ 
= &  \Sh_1\left(  \sum_{I}  \sum_{J} \pm \widehat{B}(I \times \tilde{J})  \inn{\Sh_2 \inn{f}{h_I}}{h_{J}} \cfn{I} |I|^{-1} \otimes  h_{J}
|\tilde{J}|^{-1/2}
\right) \\ 
= &  \Sh_1\left(  \sum_{I}  \sum_{J} \pm \widehat{B}(I \times \tilde{J})  \inn{\Sh_2 f}{h_I\otimes h_{J}} h_I^1  \otimes  h_{J}
|I|^{-1/2}|\tilde{J}|^{-1/2}\right). 
\end{align*}
 $\Sh_1( P^1_B(\Sh_2 f  )$.

{\bf Case $\mathbf{I \subsetneq K , J \subsetneq L}$.} Here we have
\begin{align*}
& \Sh_1 \left ( \sum_{K} \sum_{I \subsetneq K}  \sum_{L} \sum_{J \subsetneq L} \widehat{B}(I \times J) \widehat{f}(K \times L) h_I h_K \otimes  h_J \Sh_2
h_L  \right ) \\
= & \Sh_1 \left ( \sum_{J,K} \sum_{I \subsetneq K} \widehat{B}(I \times J) h_I h_K 
\otimes \left ( \sum_{L \supsetneq J} \inn{\inn{f}{h_K}}{h_L}   \Sh_2 h_L \cfn{J}\right ) h_J   \right ).
\end{align*}
By using the definition of the shift, and the known average identity $\dsp \inn{ f }{h_J^1}|J|^{-1/2} = \sum_{I \supsetneq J} \hcf{I} h_I
\cfn{J}$, we have
\begin{align*}
& \sum_{L \supsetneq J} \inn{\inn{f}{h_K}}{h_L}  \Sh_2 h_L \cfn{J}  = \Sh_2\left( \sum_{L \supsetneq J} \inn{\inn{f}{h_K}}{h_L} h_L
\right) \cfn{J} \\ 
& = \sum_{L \supseteq J} a_L \inn{\inn{f}{h_K}}{h_{\tilde{L}}} h_L \cfn{J}\\
& = a_J \inn{\inn{f}{h_K}}{h_{\tilde{J}}} h_J + \sum_{L \supsetneq J} a_L \inn{\inn{f}{h_K}}{h_{\tilde{L}}} h_L \cfn{J} \\
& = \inn{\Sh_2 \inn{f}{h_K}}{h_J^1} |J|^{-1/2} \cfn{J} + \inn{\Sh_2 \inn{f}{h_K}}{h_J}h_J .
\end{align*}
This divides the original sum into two sums $S_1 + S_2$.  The first one, $S_1$, is equal to
\begin{align*}
& \Sh_1 \left ( \sum_{K} \sum_{I \subsetneq K}  \sum_{J}  \widehat{B}(I \times J) \inn{\Sh_2 \inn{f}{h_K}}{h_J^1} h_I h_K  
\otimes \frac{h_J} { |J|^{\frac{1}{2}} } \right ) \\
= &  \Sh_1 \left ( \sum_{I}   \sum_{J}  \widehat{B}(I \times J) \left ( \sum_{K \supsetneq I}\inn{ \inn{\Sh_2 f}{h_J^1}}{h_K} h_K \cfn{I} 
\right ) h_I \otimes  \frac{h_J} { |J|^{\frac{1}{2}} }  \right ) \\
= &  \Sh_1 \left ( \sum_{I}   \sum_{J}  \widehat{B}(I \times J) \inn{ \inn{\Sh_2 f}{h_J^1}}{h_I^1}  \frac{h_I 
\otimes  h_J }{ |I|^{ \frac{1}{2} } |J|^{ \frac{1}{2} } } \right ) \\
= &  \Sh_1 \left ( \sum_{I}   \sum_{J}  \widehat{B}(I \times J) \inn{\Sh_2 f}{h_I^1 \otimes h_J^1}  \frac{h_I 
\otimes  h_J }{ |I|^{ \frac{1}{2} } |J|^{ \frac{1}{2} } }
\right ).
\end{align*}
Which has the form $\Sh_1( P^3_B(\Sh_2 f  ) )$.  And with similar computations, we get
\begin{align*}
S_2 & =  \Sh_1 \left ( \sum_{I}   \sum_{J}  \widehat{B}(I \times J) \inn{\Sh_2 f}{h_I^1 \otimes h_J}  h_I 
\otimes  h_J^1 |I|^{-1/2} |J|^{-1/2} \right ) \\
& = \Sh_1( P^5_B(\Sh_2 f  ) ).
\end{align*}

{\bf Case $\mathbf{I = K , J \subsetneq L}$.} In this case we get
\begin{align*}
& \Sh_1 \left ( \sum_{I}  \sum_{L} \sum_{J \subsetneq L} \widehat{B}(I \times J) \widehat{f}(I \times L)  h_I^2 
\otimes  h_J \Sh_2 h_L  \right ) \\
=\ & \Sh_1 \left ( \sum_{I} \sum_{J}  \widehat{B}(I \times J) h_I^2 
\otimes \left ( \sum_{L \supsetneq J} \inn{\inn{f}{h_I}}{h_L}\Sh_2 h_L \cfn{J} \right ) h_J   \right ) \\
=\ & \Sh_1 \left ( \sum_{I} \sum_{J}  \widehat{B}(I \times J) h_I^2 
\otimes \inn{\Sh_2 \inn{f}{h_I}}{h_J^1} h_J|J|^{-1/2}    \right ) \\
& + \Sh_1 \left ( \sum_{I} \sum_{J}  \widehat{B}(I \times J) h_I^2 
\otimes   \inn{\Sh_2 \inn{f}{h_I}}{h_J}h_J  \right ) \\
=\ & S_1 + S_2.
\end{align*}
Again, by the definition of the shift 
\begin{align*}
S_1 & = \Sh_1 \left ( \sum_{I} \sum_{J}  \widehat{B}(I \times J) h_I^2 
\otimes \inn{\Sh_2 \inn{f}{h_I}}{\cfn{J}|J|^{-1}}  h_J   \right ) \\
 & = \Sh_1 \left ( \sum_{I} \sum_{J}  \widehat{B}(I \times J)   \inn{\Sh_2 f } {h_I 
\otimes \cfn{J}|J|^{-1}} \cfn{I} |I|^{-1} \otimes  h_J  \right ) \\
 & = \Sh_1 \left ( \sum_{I} \sum_{J}  \widehat{B}(I \times J)   \inn{\Sh_2 f }{h_I 
\otimes h^1_J} h_I^1 \otimes  h_J  |I|^{-1/2} |J|^{-1/2} \right ).
\end{align*}
Which has the form $\Sh_1( P^4_B(\Sh_2 f  ) )$.  And similarly
\begin{align*}
S_2 & = \Sh_1 \left ( \sum_{I} \sum_{J}  \widehat{B}(I \times J)   \inn{\Sh_2 f }{h_I 
\otimes h_J} h_I^1 \otimes  h_J^1  |I|^{-1/2} |J|^{-1/2} \right ) \\
& = \Sh_1( (P^3_{B^*})^*(\Sh_2 f  ) ).
\end{align*}

{\bf Case $\mathbf{I \subsetneq K ,  J = L}$.} last case we have
\begin{align*}
& \Sh_1 \left ( \sum_{K}  \sum_{J} \sum_{I \subsetneq K} \widehat{B}(I \times J) \widehat{f}(K \times J)  h_I h_K \otimes  h_J \Sh_2 h_J  \right ) \\
& \Sh_1 \left ( \sum_{I}  \sum_{J}  \widehat{B}(I \times J) \left( \sum_{K \supsetneq I} \inn{\inn{f}{h_J}}{h_K} h_K \cfn{I} \right) h_I  \otimes  h_J
\Sh_2 h_J \right ) \\
& \Sh_1 \left ( \sum_{I}  \sum_{J}  \widehat{B}(I \times J)  \inn{\inn{f}{h_J}}{h_I^1}  h_I |I|^{-1/2}  \otimes  (h_{J_-} - h_{J_+}) |J|^{-1/2} \right ).
\end{align*}
This is a sum of two terms of the form
$$
\Sh_1 \left ( \sum_{I}  \sum_{J}  \pm \widehat{B}(I \times \tilde{J})  \inn{f}{h_I^1 \otimes h_{\tilde{J}}}  \frac{h_I  \otimes  h_J}{|I|^{\frac{1}{2}} |\tilde{J}|^{\frac{1}{2}} } \right )  = \Sh_1 (P^2_B(f) ).
$$
This concludes the proof of the estimate for the term $\widetilde{T}_2$.

\subsection{Remark: Logarithmic estimate}

Note that, because of (\ref{bmotrace}), the previous estimates for the upper bound depend on a dimensional constant.    Using a
slightly different ordering of the terms in the formal Haar expansion of the product $B f$, we obtain a decomposition in paraproducts of the form
\begin{align*}
& \ \ \ \sum_{R \in \cD^2} \inn{B}{h_R^{(0,0)}} \inn{f}{h_R^{(0,0)}}h_R^{(1,1)} \! + \! \pp{0}{0}{0}{1}{1}{0} \\
& + \pp{0}{1}{0}{0}{1}{0}  \! + \! \sum_{R \in \cD^2} \inn{B}{h_R^{(1,0)}} \inn{f}{h_R^{(0,0)}}h_R^{(0,1)} \\
& + \pp{1}{0}{0}{1}{0}{0} \! + \! \pp{1}{1}{0}{0}{0}{0} \\
&  + \sum_{R \in \cD^2} \inn{B}{h_R^{(0,0)}} \inn{f}{h_R^{(1,0)}}h_R^{(0,1)} \! + \! \pp{0}{0}{1}{1}{0}{0} \\
& + \pp{0}{1}{1}{0}{0}{0} \\
& = (T_1 +T_2 + T_3 + T_4 + T_5 + T_6 + T_7 + T_8 + T_9)(f).
\end{align*}
Here, $h_R^{(\eps, \delta)} = h_I^{\eps} h_J^{\delta}$, with $\eps, \delta \in \{ 0,1 \}$, and $h_I^0 = h_I$, $h_I^1 = |I|^{-1/2}\cfn{I}$. Then, 
$$
\com{\com{M_B}{\Sh_1}}{\Sh_2}(f) = \shc{T_1}(f) + \cdots \shc{T_9}(f).
$$
Therefore, to find an upper bound for the commutator, it suffices to find upper bounds for the different paraproducts in the above expansion.  By the
previous section, this upper bound depends also on a dimensional constant, however, it is possible for the terms $T_1$, $T_6$, and $T_8$ (by duality),
to find a better estimate of order $\log^2(1+d)$.  This is possible due to a generalization of the results obtained by Pisier in \cite{Pi2005} for the
one parameter case, combined with the characterization by two index martingales given by Bernard in \cite{Be1979}.

With the rest of the terms, it's still not clear how to find this improved dimensional bound for the paraproduct, since we don't have a representation in two-index martingales in these cases, or the appropriate embedding theorem.

\section{Lower bound}

The lower bound can be proved by using the result in the scalar case (proved by Ferguson and Lacey
in \cite{FeLa2002}).   That, is, there is a constant $C > 0$ such that
$$
\| b \|_{BMO} \leq C \| [[M_b, H_1],H_2] \|_{L^2 \en L^2},
$$
for all scalar functions $b$ in $BMO(\bR^2)$.  Let us recall the definition of $BMO$ given in \ref{def:BMO}.  The
lower bound estimate in the matrix-valued setting is the following
\begin{teo}[Lower bound]
Let $B$ be a matrix-valued function on $\bR^2$, then
$$
d^{-2}  \| B \|_{BMO} \lesssim  \left\| [[M_B, H_1],H_2] \right\|_{L^2(\bC^d) \en L^2(\bC^d)}.
$$
\end{teo}

{\bf Proof: } Denote by $\widehat{B}(R)$ the wavelet coefficient $\inn{B}{v_R}$.  Consider the functions $f, g \in L^2(\bC)$.  Let $\{ \row{\vect{e}}{1}{d} \}$
represent the cannonical basis of $\bR^d$, then, for $1 \leq i,j \leq d$, the functions $\tilde f =
f \vect{e}_i$ and $\tilde g = g \vect{e}_j$ both belong to $L^2(\bC^d)$.  If $B = (b_{ij})$, an easy
computation shows that 
$$
\inn{[[M_B, H_1], H_2]\tilde f}{\tilde g}_{L^2(\bC^d)} = \inn{[[M_{b_{ji}}, H_1], H_2] f }{ g
}_{L^2(\bC)}
$$
Therefore, for every $i,j \in \{ 1, \ldots, d \}$, we have
\begin{equation}\label{CommScalartoMatrix}
\| [[M_{b_{ji}}, H_1], H_2] \|_{L^2(\bC) \en L^2(\bC)} \leq \|  [[M_B, H_1], H_2] \|_{L^2(\bC^d)
\en L^2(\bC^d)}.
\end{equation}
Let $\{ E_{ij} : 1\leq i,j \leq d\}$ be the canonical basis for the $d \times d$ matrices, that is,
$(E_{ij})_{kl} = \delta_{ik}\delta_{jl}$.  We can write $ B = \sum_{i,j} b_{ij}E_{ij}$, and proceed
to find an estimate for the $BMO$ norm of the matrices $\tilde{B}_{ij} = b_{ij}E_{ij}$.  

Note that $\widehat{\tilde{B}}_{ij}(R) \widehat{\tilde{B}}_{ij}(R)^* = \widehat{\tilde{B}}_{ij}(R)^*
\widehat{\tilde{B}}_{ij}(R) = \widehat{b}_{ij}(R)E_{ij}  \overline{\widehat{b}}_{ij}(R)E_{ji} =
|\widehat{b}_{ij}(R)|^2 E_{ii}$.  Then, for any open set $ U \subc \bR^2$, we have
\begin{align*}
\frac{1}{|U|} \sum_{R \subc U}  \widehat{\tilde{B}}_{ij}(R) \widehat{\tilde{B}}_{ij}(R)^* & =
\frac{1}{|U|} \sum_{R \subc U}   |\widehat{b}_{ij}(R)|^2 E_{ii} \\
& \leq \frac{1}{|U|} \sum_{R \subc U}
  |\widehat{b}_{ij}(R)|^2 I_d \leq \| b_{ij} \|_{BMO} I_d.
\end{align*}
Using the one parameter result, and equation \ref{CommScalartoMatrix}, we get 
\begin{align*}
\frac{1}{|U|} \sum_{R \subc U}  \widehat{\tilde{B}}_{ij}(R) \widehat{\tilde{B}}_{ij}(R)^* & 
\lesssim \| [[M_{b_{ji}}, H_1], H_2] \|_{L^2(\bC) \en L^2(\bC)} I_d \\
& \leq \|  [[M_B, H_1], H_2] \|_{L^2(\bC^d) \en L^2(\bC^d)}.
\end{align*}
That is, $\dsp \| \tilde{B}_{ij} \|_{BMO} \lesssim \|  [[M_B, H_1], H_2]
\|_{L^2(\bC^d) \en L^2(\bC^d)}$.  Therefore,
$$
\| B \|_{BMO} \leq \sum_{i,j} \| \tilde{B}_{ij} \|_{BMO} \lesssim d^2 \| [[M_B, H_1], H_2]
\|_{L^2(\bC^d) \en L^2(\bC^d)}.
$$
Which is the desired lower bound. \hfill $\blacksquare$

\end{document}